\documentclass[a4paper]{amsart}

\usepackage[mathscr]{eucal}
\usepackage{srcltx}
\usepackage{epsfig}
\usepackage{amssymb}
\usepackage{longtable}
\usepackage{amssymb}
\usepackage{amscd,amsfonts}
\usepackage{amsthm}
\usepackage{amsmath}
\usepackage{latexsym}

\theoremstyle{plain}
\newtheorem{theorem}{Theorem}

\newtheorem{lemma}[theorem]{Lemma}

\theoremstyle{definition}
\newtheorem{definition}[theorem]{Definition}

\newtheorem{example}[theorem]{Example}

\numberwithin{theorem}{section}

\newcommand{\Ch}{\mathbf{C}}

\newcommand{\card}{\mbox{card}}

\newcommand{\Proof}{{\it Proof. \,}}

\def\Qco{\mathfrak{Qco}}

\newcommand{\OHom}[3]{\mbox{Hom}_{#1}(#2,#3)}
\newcommand{\Ext}[4]{\mbox{Ext}^{#1}_{#2}(#3,#4)}

\newcommand{\rmod}[1]{\mbox{\rm{Mod}--}{#1}}

\newcommand{\Spec}{\mbox{\rm Spec}}

\begin{document}

\title[Descent of restricted flat Mittag--Leffler modules]
{Descent of restricted flat Mittag--Leffler modules and generalized vector bundles}
\author{Sergio Estrada}
\address{Departamento de Matem\'atica Aplicada, Universidad de Murcia\\
Campus del Espinardo 30100, Spain}
\email{sestrada@um.es}
\author{Pedro Guil Asensio}
\address{Departamento de Matem\'aticas, Universidad de Murcia\\
Campus del Espinardo 30100, Spain}
\email{paguil@um.es}
\author{Jan Trlifaj}
\address{Charles University, Faculty of Mathematics and Physics, Department of
Algebra \\
Sokolovsk\'{a} 83, 186 75 Prague 8, Czech Republic}
\email{trlifaj@karlin.mff.cuni.cz}
\thanks{The first author partially supported by DGI MTM2008-03339, by the Fundaci\'on Seneca
and by the Junta de Andaluc\'{\i}a, Consejer\'{\i}a de Econom\'{\i}a, Innovaci\'on y Ciencia 
and FEDER funds. The third author supported by GA\v CR 201/09/0816 and MSM 0021620839}


\date{\today}
\begin{abstract} 
A basic question for any property of quasi--coherent sheaves on a scheme $X$ is whether the property is local, that is, it can be defined using any open affine covering of $X$. Locality follows from the descent of the corresponding module property: for (infinite dimensional) vector bundles and Drinfeld vector bundles, it was proved by Kaplansky's technique of d\'evissage already in \cite[II.\S3]{RG}. Since vector bundles coincide with $\aleph_0$--restricted Drinfeld vector bundles, a question arose in \cite{EGPT} of whether locality holds for $\kappa$--restricted Drinfeld vector bundles for each infinite cardinal $\kappa$. We give a positive answer here by replacing the d\' evissage with its recent refinement involving $\mathcal C$--filtrations and the Hill Lemma.  
\end{abstract}

\maketitle 
\section{Introduction}
In \cite[\S2]{D} Drinfeld considered generalizations of the notion of a finite dimensional vector bundle on a scheme to the infinite dimensional setting.  
One of the options suggested replaced finitely generated projective modules by the flat Mittag--Leffler ones, leading thus to what was later called a Drinfeld vector bundle in \cite{EGPT}. Mittag--Leffler modules have been studied for decades, starting from the classic works of Grothendieck, Raynaud and Gruson \cite{RG} until recent contributions in \cite{AH}, \cite{HT} et al. 

Another ingredient of the theory is due to Quillen. Following \cite{Q}, one can
compute morphisms between two objects $A$ and $B$ of the derived category of the
category $\Qco (X)$ of quasi--coherent sheaves on a scheme $X$. First, one
introduces a model category structure on $\Ch (\Qco (X))$ (= the category of
unbounded chain complexes on $\Qco (X)$). Morphisms between $A$ and $B$ are then
computed as the $\Ch (\Qco (X))$--morphisms between cofibrant and fibrant
replacements of $A$ and $B$ modulo chain homotopy. Much later, Hovey and
Gillespie \cite{H}, \cite{G} developed a machinery for transferring the
computation to $\Ch (\Qco (X))$ using complete cotorsion pairs. 

A test question for the Drinfeld options is to what extent they fit the Quillen--Hovey theory. While the answer for the general setting of arbitrary flat Mittag--Leffler modules is negative, there is a positive answer in case we admit only filtrations by Mittag--Leffler modules of restricted size \cite{EGPT}. This yields the notion of a $\kappa$--restricted Drinfeld vector bundle for each infinite cardinal $\kappa \geq \aleph_0$. In fact, the $\aleph_0$--restricted Drinfeld vector bundles are exactly the (infinite dimensional) vector bundles from \cite[\S2, Definition]{D}, that is, those quasi--coherent sheaves all of whose modules of sections are projective, but not necessarily finitely generated.     

\medskip
While defining the notion of a restricted Drinfeld vector bundle, one faces the question of whether the notion is local, i.e., it can be defined using any open affine covering of the scheme. The question is known to have a positive answer in the two extreme cases: for all Drinfeld vector bundles, and for all (infinite dimensional) vector bundles. In \cite[Seconde partie]{RG} it was proved that locality follows from the descent of the corresponding module property, that of a flat Mittag--Leffler module, and a projective module, respectively. Raynaud and Gruson proved the descent of projectivity using Kaplansky's technique of d\'evissage: a~module $M$ is projective, if and only if $M$ is both flat Mittag--Leffler and a direct sum of countably generated modules (cf.\ \cite[\S2]{D}). However, whether the notion of a $\kappa$--restricted Drinfeld vector bundle is local also for $\kappa > \aleph_0$ remained open \cite{EGPT}. Our main result here gives a positive answer:

\begin{theorem}\label{local} Let $\kappa$ be an infinite cardinal. Then the notion of a $\kappa$--restricted Drinfeld vector bundle is local. That is, if $X$ is any scheme, $V$ a set of open affine sets of $X$ such that $X = \bigcup_{v \in V} v$, and $\mathscr M$ a quasi--coherent sheaf on $X$ such that $\mathscr M (v)$ is a $\kappa$--restricted flat Mittag--Leffler $\mathscr R (v)$--module for all $v \in V$, then $\mathscr M$ is a $\kappa$--restricted Drinfeld vector bundle on $X$.
\end{theorem}

Our proof is also based on proving descent, this time for restricted flat Mittag--Leffler modules, but replaces Kaplansky's d\'evissage by a more subtle technique for dealing with filtrations due to Hill. On the way, we also notice that the $\kappa$--restricted flat Mittag--Leffler modules for $\kappa > \aleph_0$ are related to another interesting class of modules, the Whitehead ones. Namely, as shown by Eklof and Shelah \cite{ES}, there exist local PIDs $R_\kappa$ such that the $\kappa$--restricted Mittag--Leffler $R_\kappa$--modules coincide with the Whitehead $R_\kappa$--modules under the assumption of G\" odel's Axiom of Constructibility (see Example \ref{proj-white}(ii)). So in this particular case the $\kappa$--restricted Drinfeld vector bundles over the affine scheme 
$\Spec (R_\kappa)$ are exactly the Whitehead vector bundles.  

\section{Preliminaries}
Let $R$ be a ring. Following \cite{RG}, we call a (right $R$--) module $M$ \emph{Mittag--Leffler} provided that the canonical map 
$M \otimes_R \prod_{i \in I} M_i \to \prod_{i \in I} M \otimes_R M_i$ is
monic for each family of left $R$--modules $( M_i \mid i \in I )$. 

As in \cite[\S2]{D}, we will primarily be interested in flat Mittag--Leffler modules, that is, the Mittag--Leffler modules such that the functor $M \otimes_R -$ is exact.
These modules were characterized in \cite[Seconde partie, Corollary 2.2.2]{RG} as the modules $M$ such that each finite subset of $M$ is contained in a countably generated pure and projective submodule of $M$. In fact, they can also be characterized without referring to purity, as the $\aleph_1$--projective modules \cite{HT}. For example, if $R$ is a PID, then a module $M$ is flat Mittag--Leffler, iff each countably generated submodule of $M$ is free.   

The corresponding notion for quasi--coherent sheaves is as follows (see \cite{EGPT}): a quasi--coherent sheaf $\mathscr F$ on a scheme $X$ is a \emph{Drinfeld vector bundle} provided that for each open affine set $u$, the $\mathscr R (u)$--module of sections $\mathscr M (u)$ is flat and Mittag--Leffler.   

We also recall \cite[\S2, Definition]{D} that a quasi--coherent sheaf $\mathscr F$ on a scheme $X$ is an (infinite dimensional) \emph{vector bundle} provided that for each open affine set $u$, the $\mathscr R (u)$--module $\mathscr M (u)$ is projective. 

A module $M$ is a \emph{Whitehead module} provided that $\Ext 1RMR = 0$ \cite{EM}. A quasi--coherent sheaf $\mathscr F$ on a scheme $X$ is a \emph{Whitehead vector bundle} provided that for each open affine set $u$, the $\mathscr R (u)$--module $\mathscr M (u)$ is Whitehead. 

Of course, each projective module is Whitehead and flat Mittag--Leffler. The relations between the classes of all Whitehead and flat Mittag--Leffler modules depend on the structure of the ring $R$. For example, if $R$ is a Dedekind domain which is not a complete discrete valuation ring, then all Whitehead modules are flat Mittag--Leffler. But if $R$ is a complete discrete valuation domain, then Whitehead modules coincide with the flat ones, so the inclusion of the classes is reversed, cf.~\cite{EM}.
  
\medskip
In \cite{EGPT} it was shown that in general, the class of all Drinfeld vector bundles does not fit the Quillen--Hovey setting for computing cohomology via model category structures. Various restricted versions of the notion were suggested there to fix the problem. In order to define them, we recall the notion of a $\mathcal C$--filtration.

\medskip
Let $\mathcal C$ be a class of modules. A module $M$ is \emph{$\mathcal C$--filtered} provided there exist  
an ordinal $\sigma$ and a chain of submodules of $M$ 
$$0 = M_0 \subseteq \dots \subseteq M_\alpha \subseteq M_{\alpha +1} \subseteq \dots \subseteq M_\sigma = M$$
which is continuous (i.e., $M_\alpha = \bigcup_{\beta < \alpha} M_\beta$ for all limit ordinals $\alpha \leq \sigma$), and for each $\alpha < \sigma$, the module $M_{\alpha +1}/M_\alpha$ is isomorphic to an element of $\mathcal C$.

The chain $( M_\alpha \mid \alpha \leq \sigma)$ is called a \emph{$\mathcal C$--filtration} of $M$. Notice that for $\sigma = 2$, $M_2$ is just the extension of the module $M_1 \in \mathcal C$ by the module $M_2/M_1 \in \mathcal C$. Also, the class of all $\mathcal C$--filtered modules includes arbitrary direct sums of copies of modules in $\mathcal C$.   

A class of modules $\mathcal A$ is \emph{closed under transfinite extensions} provided that $M \in \mathcal A$ whenever $M$ is $\mathcal A$--filtered. 
In this case $\mathcal A$ is closed under extensions and arbitrary direct sums. 

For example, if $\mathcal A = {}^\perp \mathcal B$ for a class $\mathcal B$, then $\mathcal A$ is closed under transfinite extensions by the Eklof Lemma  \cite[3.1.2]{GT}. Here, we define  ${}^\perp \mathcal B = \mbox{Ker} \Ext 1R{-}{\mathcal B}  = \{ M \in \rmod R \mid \Ext 1RMN = 0 \mbox{ for all } N \in \mathcal B \}$. Similarly the class $\mathcal B ^\perp$ is defined; a pair of classes $(\mathcal A, \mathcal B)$ is a \emph{cotorsion pair} provided that $\mathcal A= {}^\perp \mathcal B$ and $\mathcal B = \mathcal A ^\perp$.  

In particular, the classes of all projective and flat modules are closed under transfinite extensions and direct summands, because they are the first components of cotorsion pairs. The class of all flat Mittag--Leffler modules is also closed under transfinite extensions and direct summands, but it is not of the form 
$^\perp \mathcal B$ for any class of modules $\mathcal B$ in general (see \cite{AH} and \cite{HT}). 

\medskip
Given a class $\mathcal A$ closed under transfinite extensions and an infinite cardinal $\kappa$, we let $\mathcal A ^{\leq \kappa}$ denote the class of all $\leq \kappa$--presented modules from $\mathcal A$, and $\mathcal A (\kappa)$ the class of all $\mathcal A ^{\leq \kappa}$--filtered modules.
If $M \in \mathcal A ^{\leq \kappa}$, then $M$ is called the \emph{$\kappa$--restricted module in $\mathcal A$}. We have the chain
 
$$\mathcal A (\aleph_0) \subseteq \mathcal A (\aleph_1) \subseteq \dots \subseteq \mathcal A (\kappa) \subseteq \mathcal A (\kappa^+) \subseteq \dots \subseteq \mathcal A = \bigcup_{\aleph_0 \leq \kappa} \mathcal A (\kappa).$$  

\begin{example}\label{threecases} 1. If $\mathcal A$ is the class of all projective modules over any ring $R$, then the chain above is constant, because  $\mathcal A (\aleph_0) = \mathcal A$ by a classic result of Kaplansky \cite{K}.

2. If $\mathcal A$ is the class of all flat modules over a ring $R$, then the chain stabilizes at $\card R + \aleph_0$, that is, $\mathcal A (\kappa) = \mathcal A$ for each infinite cardinal $\kappa \geq \card R$. This follows from Enochs' solution of the Flat Cover Conjecture \cite{BEE}.

3. Let $\mathcal A$ be the class of all flat Mittag--Leffler modules over any ring $R$. If $\kappa$ is an infinite cardinal, then the modules in the class $\mathcal A (\kappa)$ are the $\kappa$--restricted flat Mittag--Leffler modules. A quasi--coherent sheaf $\mathscr F$ on a scheme $X$ is called \emph{$\kappa$--restricted Drinfeld vector bundle} provided that for each open affine set $u$, the $\mathscr R (u)$--module of sections $\mathscr M (u)$ is $\kappa$--restricted flat Mittag--Leffler. 

If $R$ is right perfect, then $\mathcal A$ is just the class of all projective modules, so we are in case 1. But if $R$ is a non--right perfect ring, then the chain above is strictly increasing starting from $\card R + \aleph_0$, because for each infinite cardinal $\kappa \geq \card R$, there exists a module $M_{\kappa^+} \in \mathcal A ^{\leq \kappa^+} \setminus  \mathcal A (\kappa)$ (see \cite[Theorem 6.10]{HT}).
\end{example}

Returning to the general case, we observe that the classes $\mathcal A (\kappa)$ are often the first components of complete cotorsion pairs. Here, a cotorsion pair $(\mathcal A, \mathcal B)$ is said to be \emph{complete} provided that for each module $M$ there is a short exact sequence $0 \to M \to B \to A \to 0$. 

\begin{lemma}\label{basic} Let $R$ be a ring and $\kappa$ be an infinite cardinal. Let $\mathcal A$ be a class of modules closed under transfinite extensions and direct summands, and such that $R \in \mathcal A$. Then $(\mathcal A (\kappa), (\mathcal A ^{\leq \kappa})^\perp)$ is a complete cotorsion pair. 
\end{lemma}
\noindent\Proof First, we observe that $\mathcal A (\kappa)^\perp = (\mathcal A ^{\leq \kappa})^\perp$ by the Eklof Lemma \cite[3.1.2]{GT}. Let 
$\mathcal E = {}^\perp ((\mathcal A ^{\leq \kappa})^\perp)$. Since $\mathcal A ^{\leq \kappa}$ has a representative set $\mathcal S$ of objects and contains the regular module $R$, \cite[3.2.4]{GT} implies that $\mathcal E$ consists of all direct summands of elements of $\mathcal A (\kappa)$. In particular, $\mathcal E \subseteq \mathcal A$ by assumption. Moreover, the Kaplansky theorem for cotorsion pairs \cite[4.2.11]{GT} shows that each element of $\mathcal E$ is 
$\mathcal E ^{\leq \kappa}$--filtered, hence $\mathcal A ^{\leq \kappa}$--filtered, and $\mathcal E = \mathcal A (\kappa)$. Finally, since $(\mathcal A ^{\leq \kappa})^\perp = \mathcal S ^\perp$, the cotorsion pair $(\mathcal A (\kappa), (\mathcal A ^{\leq \kappa})^\perp)$ is complete by  \cite[3.2.1]{GT}.\qed   

\medskip
Moreover, if the class $\mathcal A ^{\leq \kappa}$ contains all its $\leq \kappa$--presented syzygies, then the cotorsion pair in Lemma \ref{basic} is hereditary, that is, $\Ext iRAB = 0$ for all $i \geq 2$, $A \in \mathcal A (\kappa)$, and $B \in (\mathcal A ^{\leq \kappa})^\perp$. This is the case when

(i) $\mathcal A$ = the class of all flat Mittag--Leffler modules over any ring by \cite{AH}, and when 

(ii) $\mathcal A$ = the class of all flat modules over any right $\kappa$--noetherian ring.

\medskip
Finally, let us consider the relations between restricted Mittag--Leffler modules and the better known classes of projective and Whitehead modules: 

\begin{example}\label{proj-white} 
(i) By the characterization of flat Mittag--Leffler modules in \cite{RG} mentioned above, all countably generated flat Mittag--Leffler modules are projective. So the notion of an $\aleph_0$--restricted flat Mittag--Leffler module coincides with that of a projective module, and (infinite dimensional) vector bundles are exactly the $\aleph_0$--restricted Drinfeld vector bundles. 

(ii) Let $\kappa$ be an infinite cardinal. By \cite[Theorem 6]{ES} there exists a local, but not complete, PID $R = R_\kappa$ of cardinality $2^\kappa$ such that for each $\leq \kappa$--generated module $M$, $M$ is Whitehead, iff $M$ is flat Mittag--Leffler. If we assume V = L (G\" odel's Axiom of Constructibility), 
then for each module $M$, we have that $M$ is Whitehead, iff $M$ is $\kappa$--restricted flat Mittag--Leffler (see e.g. \cite[Theorem 10.1.5]{GT}). So in this case the $\kappa$--restricted Drinfeld vector bundles over the affine scheme $X = \Spec (R_\kappa)$ are exactly the Whitehead vector bundles over $X$.             
\end{example}

\section{Locality for induced quasi--coherent sheaves}
Let $\mathfrak P$ be a property of modules. For each commutative ring $R$, we let $\mathfrak P _R$ denote the class of all $R$--modules satisfying $\mathfrak P$. Throughout, we will assume that $\mathfrak P$ is compatible with ring direct products in the following sense: if $R = \prod_{i < n} R_i$ and $M_i \in \mathfrak P _{R_i}$ for all $i < n$, then $M = \prod_{i<n} M_i \in \mathfrak P _R$. 

\begin{definition}\label{induced}
For a scheme $X$, we define the quasi--coherent sheaf \emph{induced by $\mathfrak P$} on $X$ (or a \emph{$\mathfrak P$--quasi--coherent sheaf on $X$}) as the quasi--coherent sheaf $\mathscr M$ such that for each open affine set $u$ of $X$, the $\mathscr R (u)$--module of sections $\mathscr M (u)$ satisfies $\mathfrak P$ (that is, $\mathscr M (u) \in \mathfrak P _{\mathscr R (u)}$).
\end{definition}  

\begin{example}\label{allex} Let $\mathfrak P$ be the property of being a projective, flat, flat Mittag--Leffler, and $\kappa$--restricted flat Mittag--Leffler module (where $\kappa$ is an infinite cardinal). Then $\mathfrak P _R$ is the class of all projective, flat, flat Mittag--Leffler, and $\kappa$--restricted flat Mittag--Leffler $R$--modules (and $\mathfrak P$ is obviously compatible with ring direct products). Moreover, $\mathfrak P$--quasi--coherent sheaves on a scheme $X$ are exactly the (infinite dimensional) vector bundles, flat quasi--coherent sheaves, Drinfeld vector bundles, and $\kappa$--restricted Drinfeld vector bundles on $X$, respectively. 
\end{example}

\begin{definition}\label{localprop} The notion of a $\mathfrak P$--quasi--coherent sheaf is \emph{local} in case for each open affine covering $X = \bigcup_{v \in V} v$ of $X$ and each quasi--coherent sheaf $\mathscr M$ on $X$, $\mathscr M$ is $\mathfrak P$--quasi--coherent provided that 
$\mathscr M (v) \in \mathfrak P _{\mathscr R (v)}$ for all $v \in V$. 
\end{definition}

So locality means that the property of being a $\mathfrak P$--quasi--coherent sheaf can be tested on any open affine covering of $X$.     

We present a useful classic tool for proving locality based on ascent and descent of the coresponding module property (cf.\ \cite{SP}): 

\begin{lemma}\label{faithflat} Let $\mathcal R$ be a class of commutative rings. Let $\mathfrak P$ be a property of modules such that for each flat homomorphism $\varphi : R \to S$ of rings in $\mathcal R$ and each $R$--module $M$, the following two conditions hold:

\begin{enumerate}
\item $M \in \mathfrak P _R$ implies $M \otimes_R S \in \mathfrak P _S$, and
\item $M \otimes_R S \in \mathfrak P _S$ implies $M \in \mathfrak P _R$ provided that $\varphi$ is faithfully flat.
\end{enumerate}

Then the notion of a $\mathfrak P$--quasi--coherent sheaf is local.
\end{lemma}
\noindent\Proof Let $X = \bigcup_{v \in V} v$ be an open affine covering of $X$ such that $\mathscr M (v) \in \mathfrak P _{\mathscr R (v)}$ for all $v \in V$. Let $u$ be an arbitrary open affine set of $X$.  Then there exists a standard open covering $u = \bigcup_{j < n} u_j$ such that for each $j < n$ there exists $v_j \in V$ such that $u_j$ is a standard open in $v_j$, so $u_j = D(f_j) = \Spec (\mathscr R (v_j)_{f_j})$ for some $f_j \in \mathscr R (v_j)$. In particular, $\mathscr M (u_j) \cong \mathscr M (v_j) \otimes_{\mathscr R (v_j)} \mathscr R (v_j) _{f_j}$, hence $\mathscr M (u_j) \in \mathfrak P _{\mathscr R (u_j)}$ by Condition (1). The compatibility of $\mathfrak P$ with ring direct products gives $\prod_{j < n} \mathscr M (u_j) \in \mathfrak P _{\prod_{j < n} \mathscr R (u_j)}$. Since the canonical morphism $\psi : \mathscr R (u) \to \prod_{j < n} \mathscr R (u_j)$ is faithfully flat, we conclude that $\mathscr M (u) \in \mathfrak P _{\mathscr R (u)}$ by Condition (2).\qed   

\begin{definition}\label{ad-property}
A property $\mathfrak P$ satisfying Conditions (1) and (2) above is said to \emph{ascend} and \emph{descend} in $\mathcal R$, respectively (see \cite{RG}). 

A property $\mathfrak P$ satisfying both Conditions (1) and (2) of Lemma \ref{faithflat} is called the \emph{ascent--descent property}, or \emph{AD--property}, in $\mathcal R$.  
\end{definition}

While the ascent is usually easy to prove, proving the descent is more involved. 

However, it is easy to see that being a flat module is an AD--property. Also the property of being a flat Mittag--Leffler module is known to be an AD-property: Condition (1) is easy while Condition (2) follows from the next lemma:

\begin{lemma}\label{MLdescent} Let $\varphi : R \to S$ be a faithfully flat homomorphism of commutative rings, and $M$ be an $R$--module such that 
$M \otimes_R S$ is a Mittag--Leffler $S$--module. Then $M$ is a Mittag--Leffler $R$--module.
\end{lemma}
\noindent\Proof First, we note that the property of being a Mittag--Leffler module can be restated in terms of the Mittag--Leffler condition for certain inverse systems: expressing $M$ as a direct limit of a direct system $( F_i, f_{ij} \mid i \leq j \in I )$ of finitely presented modules, $M = \varinjlim_{i \in I} F_i$, and applying the functor $\OHom R{-}R$, we obtain the inverse system $( \OHom R{F_i}R, \OHom R{f_{ij}}R \mid i \leq j \in I )$. Then $M$ is a Mittag--Leffler $R$--module, iff this inverse system satisfies the Mittag--Leffler condition. 

We have  $M \otimes_R S = \varinjlim_{i \in I} F_i \otimes_R S$, and $M \otimes_R S$ is a Mittag--Leffler $S$--module by assumption. The Mittag--Leffler condition says that for each $i \in I$ the family $( \mbox{Im} \OHom S{f_{ij}\otimes_R S}S \mid i \leq j \in I )$ of $S$--submodules of $\OHom S{F_i \otimes_R S}S$ stabilizes starting from some $j \geq i$. Since $F_i$ is finitely presented, there is a natural isomorphism 
$$\OHom S{F_i \otimes_R S}S \cong \OHom R{F_i}R \otimes_R S.$$ 
Since $\varphi$ is faithfully flat, for each $i \in I$ the family $( \mbox{Im} \OHom R{f_{ij}}R \mid i \leq j \in I )$ of $R$--submodules of 
$\OHom R{F_i}R$ stabilizes starting from some $j \geq i$. But the latter just says that $M$ is a Mittag--Leffler $R$--module.\qed

\section{Restricting ad--properties of modules}
In order to treat the restricted Drinfeld vector bundles, we will need to transfer the unrestricted version of Lemma \ref{faithflat} to a restricted one. The following result (known as the Hill Lemma, and proved e.g.\ in \cite[Theorem 4.2.6]{GT}) will be our tool for refining the technique of d\'evissage from $\mathcal C$--direct sums of modules to the more general setting of $\mathcal C$--filtered modules: 

\begin{lemma}\label{hill} 
Let $R$ be a ring, $\lambda$ a regular infinite cardinal, and
$\mathcal C$ a class of ${<} \lambda$--presented modules. Let $P$ be
a module with a $\mathcal C$--filtration $\mathcal P = (P_\alpha
\mid \alpha \leq \sigma )$. Then there is a family $\mathcal F$ 
consisting of submodules of $P$ such that
\begin{itemize}
\item[{\rm (H1)}] $\mathcal P \subseteq \mathcal F$,
\item[{\rm (H2)}] $\mathcal F$ is closed under arbitrary sums and intersections,
\item[{\rm (H3)}] $N^\prime/N$ is $\mathcal P$--filtered for all $N, N^\prime \in \mathcal F$ such
that $N \subseteq N^\prime$, and
\item[{\rm (H4)}] If $N \in \mathcal F$ and $T$ is a subset of $P$ of cardinality ${<} \lambda$,
then there exists $N^\prime \in \mathcal F$ such that $N \cup T \subseteq N^\prime$ and $N^\prime/N$
is ${<} \lambda$--presented.
\end{itemize}
\end{lemma}

For the restricted version, we will use the following   

\begin{theorem}\label{restricted} Let $\mathcal R$ be a class of commutative rings and $\mathfrak P$ be an AD--property in $\mathcal R$. 
Assume moreover that for each $R \in \mathcal R$, the class $\mathfrak P _R$ is closed under pure submodules and transfinite extensions, and contains $R$.  

Let $\kappa$ be an infinite cardinal such that for each $R \in \mathcal R$, each $M \in \mathfrak P _R$, and each subset $A \subseteq M$ of cardinality $\leq \kappa$, there exists a $\leq \kappa$--presented pure submodule $N \subseteq M$ containing $A$.

Let $\mathfrak R$ be the restricted property of modules defined by $\mathfrak R
_R = \mathfrak P _R (\kappa)$ for each $R\in \mathcal R$. Then $\mathfrak R$ is
an AD--property in $\mathcal R$. 
\end{theorem}
\noindent\Proof First, we prove that Condition (1) of Lemma \ref{faithflat} holds for $\mathfrak R$. If $M \in \mathfrak P _R (\kappa)$, then $M$ is 
$(\mathfrak P _R)^{\leq \kappa}$--filtered. Since $S$ is a flat $R$--module, $M \otimes_R S$ is $\mathcal C$--filtered where 
$\mathcal C$ is the class of all modules of the form $N \otimes_R S$ where $N \in (\mathfrak P _R)^{\leq \kappa}$. 
By Condition (1) for $\mathfrak P$, we have $N \otimes_R S \in \mathfrak P _S$. Since $N$ is $\leq \kappa$--presented as $R$--module, 
so is $N \otimes_R S$ as $S$--module. Thus $M \otimes_R S$ is $(\mathfrak P _S)^{\leq \kappa}$--filtered, and Condition (1) holds for $\mathfrak R$.      

In order to prove Condition (2) for $\mathfrak R$, let $M$ be an $R$--module such that $P = M \otimes_R S \in \mathfrak P _S (\kappa)$,
so there exists a $(\mathfrak P _S)^{\leq \kappa}$--filtration $\mathcal P = (P_\alpha \mid \alpha \leq \sigma )$ of $P$. By Condition (2) for $\mathfrak P$, we have $M \in \mathfrak P _R$. 

Consider the family $\mathcal F$ corresponding to the filtration $\mathcal P$ by Lemma \ref{hill} (for the regular infinite cardinal $\lambda = \kappa^+$). Notice that Condition (H3) of Lemma \ref{hill} yields for each $N \in \mathcal F$ that $N, P/N \in \mathfrak P _S (\kappa)$. 

We will use the family $\mathcal F$ to construct a $(\mathfrak P _R)^{\leq \kappa}$--filtration $(M_\beta \mid \beta \leq \tau )$ of $M$ by induction on $\beta$ as follows: 

We let $M_0 = 0$, and assume that $M_\beta$ is constructed so that $M_\beta$ is a pure submodule of $M$ such that $N = M_\beta \otimes _R S \in \mathcal F$ and $M/M_\beta \in \mathfrak P _R$. 

Assume there exists $x \in M \setminus M_\beta$. Since $M/M_\beta \in \mathfrak P _R$, the assumption on $\mathfrak P _R$ yields a $\leq \kappa$--presented submodule $U_1 = Q_1/M_\beta$ of $M/M_\beta$ containing $x + M_\beta$ such that $U_1$ is pure in $M/M_\beta$ (and hence $Q_1$ is pure in $M$ by \cite[Lemma 1.2.17(b)]{GT}). Then $N = N_0 \subseteq Q_1 \otimes_R S$.  

Condition (H4) of Lemma \ref{hill} provides us with $N_1 \in \mathcal F$ such that $Q_1 \otimes_R S \subseteq N_1$ and $N_1/N_0$ is $\leq \kappa$--presented. Again, the assumption on $\mathfrak P _R$ gives a $\leq \kappa$--presented submodule $U_1 \subseteq U_2 = Q_2/M_\beta$ of $M$ such that $U_2$ is pure in $M/M_\beta$ (and hence $Q_2$ is pure in $M$), and $N_1 \subseteq Q_2 \otimes_R S$. 
Proceeding similarly, we obtain a sequence of pure $R$--submodules $M_\beta = Q_0 \subseteq Q_1 \subseteq Q_2 \subseteq \dots$ of $M$, and a sequence $N = N_0 \subseteq N_1 \subseteq N_2 \subseteq \dots$ of elements of $\mathcal F$. 

Let $M_{\beta +1} = \bigcup_{i < \omega} Q_i$ and $N^\prime = \bigcup_{i < \omega} N_i$. Then $M_{\beta +1}$ is a pure submodule of $M$, and $N^\prime \in \mathcal F$ by Condition (H2) of Lemma \ref{hill}. Moreover, $M_{\beta + 1}/M_\beta = \bigcup_{i < \omega} U_i$ is a $\leq \kappa$--presented pure submodule of $M/M_\beta$, hence $M_{\beta + 1}/M_\beta \in (\mathfrak P _R)^{\leq \kappa}$ by the assumption on $\mathfrak P _R$. Moreover, $x \in M_{\beta + 1}$ and $M_{\beta+1} \otimes _R S = N^\prime$ by construction. 

Since $(M/M_{\beta +1}) \otimes_R S \cong P/N^\prime \in \mathfrak P _S (\kappa)$, Condition (2) for $\mathfrak P$ gives 
$M/M_{\beta +1} \in \mathfrak P _R$. 

If $\beta$ is a limit ordinal, we let $M_\beta = \bigcup_{\gamma < \beta} M_\gamma$. Then $M_\beta$ is pure in $P$, and $M_\beta \otimes_R S = \bigcup_{\gamma < \beta} (M_\gamma \otimes_R S) \in \mathcal F$ by Condition (H2) of Lemma \ref{hill}. Again, $(M/M_{\beta}) \otimes _R S \cong P/(M_\beta \otimes_R S) \in \mathcal F$, so Condition (2) for $\mathfrak P$ implies that $M/M_\beta \in \mathfrak P _R$.  

By construction, there exists an ordinal $\tau$ such that $M_\tau = M$, hence $(M_\beta \mid \beta \leq \tau )$ is the desired 
$(\mathfrak P _R)^{\leq \kappa}$--filtration of $M$.\qed   

\medskip
To complete the picture, we recall several known properties of flat Mittag--Leffler modules: 

\begin{lemma}\label{e-g-p-t} Let $R$ be a ring. Then the class of all flat
Mittag--Leffler modules is closed under pure submodules. 

Let $\kappa$ be an infinite cardinal. Then each $\leq \kappa$--generated
flat Mittag--Leffler module is $\leq \kappa$--presented. Moreover, each subset
of cardinality $\leq \kappa$ in a flat Mittag--Leffler module $M$ is contained
in a pure and $\leq \kappa$--presented submodule $N$ of $M$. 
\end{lemma}
\noindent\Proof This follows by \cite[Proposition 3.8]{BH} and \cite[Lemma 2.7(ii)]{EGPT}.\qed

\medskip
Now, we can prove our main result:

\medskip
\noindent\emph{Proof of Theorem \ref{local}.} Let $\mathfrak P$ denote the
property of being a flat Mittag--Leffler module. By Lemmas \ref{MLdescent} and
\ref{e-g-p-t}, $\mathfrak P$ satisfies the assumptions of Theorem
\ref{restricted} for $\mathcal R$ = the class of all commutative rings. Thus the
property $\mathfrak R$ of being a $\kappa$--restricted flat Mittag--Leffler
module satisfies Conditions (1) and (2) of Lemma \ref{faithflat}, and 
the claim follows.\qed

\end{document}